# Series solutions of Heun-type equation in terms of orthogonal polynomials

A. D. Alhaidari

*Saudi Center for Theoretical Physics, P. O. Box 32741, Jeddah 21438, Saudi Arabia*

**Abstract**: We introduce a nine-parameter Heun-type differential equation and obtain three classes of its solutions as series of square integrable functions written in terms of the Jacobi polynomial. The expansion coefficients of the series satisfy three-term recursion relations, which are solved in terms of orthogonal polynomials with continuous and/or discrete spectra. Some of these are well-known polynomials while others are either new or modified versions of known ones.

## I. Introduction

Recently, we introduced the following hypergeometric-type and confluent-hypergeometric-type differential equations

$$\left\{x(1-x)\frac{d^2}{dx^2}+\left[a-x(a+b)\right]\frac{d}{dx}+\frac{A_+}{x}+\frac{A_-}{1-x}+A_1x\right\}y(x)=A_0y(x), \qquad (1a)$$

$$\left[x\frac{d^2}{dx^2}+(a+bx)\frac{d}{dx}+\frac{A_+}{x}+A_1x\right]y(x)=A_0y(x), \qquad (1b)$$

where $\{a,b,A_0,A_1,A_\pm\}$ are real parameters. We obtained their solutions as infinite series in terms of orthogonal polynomials with continuous and/or discrete spectra [1]. Most of those are well-known polynomials while few are not. In this paper, we extend that work and consider the following nine-parameter Heun-type differential equation

$$\left\{\frac{d^2}{dx^2}+\left(\frac{a}{x}+\frac{b}{x-1}+\frac{c}{x-d}\right)\frac{d}{dx}+\frac{1}{x(x-1)(x-d)}\left[\frac{A}{x}+\frac{B}{x-1}+\frac{C}{x-d}+xD-E\right]\right\}y(x)=0, \qquad (2)$$

where $\{a,b,c,d,A,B,C,D,E\}$ are real parameters with $d\neq 0,1$. The original Heun equation corresponds to $A=B=C=0$ and $D=\alpha\beta$ with $\alpha+\beta+1=a+b+c$ [2-4]. A short review of the Heun equation and its use in physics and mathematics could be found in [5,6]. Our study of this problem is motivated by the work of F. A. Grünbaum *et al* on the Heun equation [7]. The Authors used the tridiagonalization approach to handle the equation. The algebraic underpinning of this approach was investigated by V. X. Genest *et al* in Ref. [8]. M. E. H. Ismail and E. Koelink used tridiagonalization in the J-matrix method to give the characteristics of the associated physical operators [9]. In this paper, we also employ a similar technique called the Tridiagonal Representation Approach (TRA) [10] to handle Eq. (2) as we did for Eq. (1a) and Eq. (1b) in [1].

The differential equation (2) has four regular singularities at $x=\{0,1,d,\infty\}$. Thus, due to Fuchs' theorem and using Frobenius method, we can write the following power series solution

$$y(x)=x^{\alpha}(x-1)^{\beta}(x-d)^{\gamma}\sum_{n=0}^{\infty}g_n(x-\hat{x})^n\,, \tag{3}$$

where $\{g_n\}$ are expansion coefficients and $\hat{x}$ is an ordinary point or one of the singular points. The exponents $\{\alpha,\beta,\gamma\}$ are real dimensionless parameters to be determined by the "indicial equation". However, we can replace $x^n$ in (3) by any polynomial $P_n(x)$ of order *n* in *x* since we can always write $x^n=\sum_{m=0}^{n}r_mP_m(x)$ for a certain set of coefficients $\{r_m\}_{m=0}^{n}$ with $r_n\neq 0$. The range of the variable *x* and singularities of Eq. (2) suggest that a proper choice would be the Jacobi polynomial $P_n^{(\mu,\nu)}(x)$ where we can rewrite (3) as follows

$$y(x)=\sum_{n=0}^{\infty}f_n\phi_n(x)\,, \tag{4}$$

with $\{f_n\}$ being new expansion coefficients and $\{\phi_n(x)\}$ is a complete set of functions whose elements are

$$\phi_n(x)=c_nx^{\alpha}(x-1)^{\beta}(x-d)^{\gamma}P_n^{(\mu,\nu)}(x)\,. \tag{5}$$

The Jacobi polynomial used here is defined in Appendix A by Eq. (A1) and the normalization constant is suggested by the orthogonality relation (A5) as $c_n=\sqrt{(2n+\mu+\nu+1)\frac{\Gamma(n+1)\Gamma(n+\mu+\nu+1)}{\Gamma(n+\mu+1)\Gamma(n+\nu+1)}}$. The real dimensionless parameters $\{\alpha,\beta,\gamma,\mu,\nu\}$ will be related to the differential equation parameters by the requirement of square integrability of $\{\phi_n(x)\}$ and by a special constraint on the expansion coefficients $\{f_n\}$ that will become clear shortly below. Square integrability is with respect to some positive measure, $d\zeta(x)$. That is, $\langle\phi_n|\phi_m\rangle\equiv\int\phi_n(x)\phi_m(x)d\zeta(x)=\delta_{nm}$. The differential equation and differential property of the Jacobi polynomial together with its recursion relation make it possible (in principle) for the action of the differential operator (2) on the basis element $\phi_n(x)$ to produce terms proportional to $\phi_n(x)$ and $\phi_{n\pm1}(x)$ with constant coefficients. For that to happen, certain relationships among the parameters of the bases and those of the differential equation and constraints thereof must be satisfied. If so, then we can write $\mathcal{J}\phi_n(x)=a_n\phi_n(x)+t_n\phi_{n+1}(x)+t_{n-1}\phi_{n-1}(x)$ where $\mathcal{J}=x(x-1)(x-d)\mathcal{D}$ and the Heun-type equation (2) is written as $\mathcal{D}y(x)=0$. Therefore, substituting *y*(*x*) given by Eq. (4) in the differential equation gives

$$\sum_{n=0}^{\infty}f_n(z)\left[a_n\phi_n(x)+t_n\phi_{n+1}(x)+t_{n-1}\phi_{n-1}(x)\right]=0\,. \tag{6}$$

Consequently, the expansion coefficients $\{f_n\}$ will satisfy the following symmetric three-term recursion relation (S3TRR) allowing us to take full advantage of the many powerful tools and analytic properties of orthogonal polynomials in the solution [11-13]

$$z\,f_n(z)=s_n\,f_n(z)+t_{n-1}\,f_{n-1}(z)+t_n\,f_{n+1}(z)\,, \tag{7}$$

for $n=1,2,3,...$ and where we wrote $a_n = s_n - z$ with $z$ being some proper function of the differential equation parameters. The recursion coefficients $\{s_n, t_n\}$ depend on these parameters and on *n* but are independent of *z* and such that $t_n^2 > 0$ for all *n*. Therefore, the solution $f_n(z)$ of (7) becomes a polynomial of degree *n* in *z* modulo an overall factor that depends on *z* but is independent of *n*. That is, we write $f_n(z) = p(z)P_n(z)$, where $p(z) = f_0(z)$ making $P_0(z)=1$. Defining $t_{-1}=0$ in (7) for $n=0$ gives $P_1(z) = (z-s_0)/t_0$ and gives the polynomials $\{P_n(z)\}$ explicitly for any degree using the S3TRR (7) and starting with the initial values $P_0(z)$ and $P_1(z)$. Thus, the solution of equation (2) becomes equivalent to the solution of the S3TRR (7).

In Appendix B, we show that if the differential equation parameters *A* and *B* are below certain critical values, $4A \le d(1-a)^2$ and $4B \le (1-d)(1-b)^2$, then we can derive three types of series solutions subject to the following constraints:

1. **Restricted solution**: corresponding to $D=0$, $4C = cd(1-d)(2-c)$, $E = \frac{A}{d} - \frac{B}{1-d} - \frac{c}{2}\left[ d\left(a+b+c-2\right) + 1 - a - \frac{c}{2} \right]$.
2. **Generalized solution**: corresponding to $4C \ge -d(1-d)(1-c)^2$.
3. **Special solution**: corresponding to $4C = cd(1-d)(2-c)$.

This is accomplished by using the TRA to derive four sets of S3TRR's associated with the Heun-type differential equation (2). In section II, we identify the orthogonal polynomial associated with the S3TTR (B10c) that corresponds to the "*restricted solution*" and show that it is the Wilson polynomial with discrete spectrum. However, the "*generalized solution*" will be derived in section III by showing that the S3TRR (B10b) is satisfied by a modified version of the Wilson polynomial with continuous and/or discrete spectrum. On the other hand, the "*special solution*" is obtained in section IV by solving the S3TRR (B10a) in terms of a new class of orthogonal polynomials whose properties (weight function, asymptotics, generating function, etc.) are yet to be derived. We ignore the S3TRR (B10d) since its solution is easily obtained from that associated with the S3TRR (B10c) by using the simple parameter map (B11).

## II. The restricted solution

For this class of solutions, the differential equation parameters $\{A,B,C,D,E\}$ must satisfy the constraints in the column headed "Case (c)" in Table 1 whereas the basis parameters $\{\alpha,\beta,\gamma,\mu,\nu\}$ are taken as shown in the first column of Table 2. We start by comparing (B10c) to the S3TRR of the normalized version of the Wilson polynomial $W_n(z^2;a,b,c,d)$ that reads as follows (see Eq. A7 in [14] and Eq. 9.1.4 in [15])

$$z^2 W_n = \left[ \frac{(n+a+b)(n+a+c)(n+a+d)(n+a+b+c+d-1)}{(2n+a+b+c+d)(2n+a+b+c+d-1)} + \frac{n(n+b+c-1)(n+b+d-1)(n+c+d-1)}{(2n+a+b+c+d-1)(2n+a+b+c+d-2)} - a^2 \right] W_n$$
$$-\frac{1}{2n+a+b+c+d-2}\sqrt{\frac{n(n+a+b-1)(n+c+d-1)(n+a+c-1)(n+a+d-1)(n+b+c-1)(n+b+d-1)(n+a+b+c+d-2)}{(2n+a+b+c+d-3)(2n+a+b+c+d-1)}}\, W_{n-1} \quad (8)$$
$$-\frac{1}{2n+a+b+c+d}\sqrt{\frac{(n+1)(n+a+b)(n+c+d)(n+a+c)(n+a+d)(n+b+c)(n+b+d)(n+a+b+c+d-1)}{(2n+a+b+c+d-1)(2n+a+b+c+d+1)}}\, W_{n+1}$$

The normalized version of the Wilson polynomial is written as (see Eq. A6 in [14] and Eq. 9.1.1 in [15])

$$W_n(z^2;a,b,c,d)= \sqrt{\left(\frac{2n+a+b+c+d-1}{n+a+b+c+d-1}\right)\frac{(a+b)_n(a+c)_n(a+d)_n(a+b+c+d)_n}{(b+c)_n(b+d)_n(c+d)_n n!}}\ {}_4F_3\left(\begin{matrix}-n,n+a+b+c+d-1,a+iz,a-iz\\ a+b,a+c,a+d\end{matrix}\middle|1\right) \tag{9}$$

where ${}_4F_3\left(\begin{matrix}a,b,c,d\\ e,f,g\end{matrix}\middle|z\right)=\sum_{n=0}^{\infty}\frac{(a)_n(b)_n(c)_n(d)_n}{(e)_n(f)_n(g)_n}\frac{z^n}{n!}$ and $(a)_n=a(a+1)(a+2)...(a+n-1)=\frac{\Gamma(n+a)}{\Gamma(a)}$. It is a polynomial of degree *n* in $z^2$. If we reparametrize the polynomial as $W_n(z^2;\sigma-\tau,\sigma+\tau,\gamma,\gamma)$, then comparing the resulting S3TRR (8) to (B10c) gives $f_n(z)=W_n(z^2;\sigma-\tau,\sigma+\tau,\gamma,\gamma)$ with the following parameters assignments

$$2\sigma=\nu+1,\quad 2\gamma=\mu+1,\quad z^2=-\tfrac{1}{4}(\mu+1)^2,\quad 2\tau=|a+b+c-1|. \tag{10}$$

For $\sigma^2\le\tau^2$, the spectrum is discrete and of finite size equal to *N*, which is the largest integer less than or equal to $|\tau|-\sigma$. The spectrum formula of the Wilson polynomial gives $z_k^2=-\left(k+\sigma-|\tau|\right)^2$, where $k=0,1,..,N$. Parameter substitution in the spectrum formula shows that the restricted solution corresponds to the discrete parameter values

$$\frac{4B_k}{1-d}=(1-b)^2-\left(2k+1-|a+b+c-1|+\sqrt{(1-a)^2-4A/d}\right)^2. \tag{11}$$

In that case, the solution for a given *k* is written as

$$y_k(x)=q(k)\sum_{n=0}^{N}W_n(z_k^2;\sigma-\tau,\sigma+\tau,\gamma,\gamma)\phi_n(x). \tag{12}$$

where $q^2(k)$ is the discrete part of the weight function of the Wilson polynomial (see Eq. C4 in [16] and Eq. 9.1.3 in [15]):

$$q^2(k)=-2\frac{\Gamma(2\sigma+2\gamma)\Gamma(2\tau)\left[\Gamma(\gamma-\sigma+\tau)\right]^2}{\Gamma(2\tau-2\sigma+1)\Gamma(2\gamma)\left[\Gamma(\gamma+\sigma+\tau)\right]^2}\left\{(k+\sigma-\tau)\frac{(2\sigma-2\tau)_k(2\sigma)_k\left[(\sigma-\tau+\gamma)_k\right]^2}{(1-2\tau)_k\left[(\sigma-\tau-\gamma+1)_k\right]^2k!}\right\}, \tag{13}$$

Writing $x(x-1)(x-d)\mathcal{D}=H-E$, we can recast the Heun-type equation, $\mathcal{D}y(x)=0$, as the following eigenvalue equation

$$H\left|y_k(x)\right\rangle=E_k\left|y_k(x)\right\rangle, \tag{14}$$

where the eigenvalues are

$$4E_k=\frac{4A}{d}-\frac{c}{2}\left[d(a+b+c-2)+1-a-\frac{c}{2}\right]-(1-b)^2 +\left(2k+1-|a+b+c-1|+\sqrt{(1-a)^2-4A/d}\right)^2 \tag{15}$$

# III. The generalized solution

For this class of solutions, the differential equation parameters $\{A,B,C\}$ and the basis parameters $\{\alpha,\beta,\gamma,\mu,\nu\}$ are to be chosen and should satisfy the constraints given in the column headed "Case (b)" in Table 1.

We start by introducing a version of the Wilson polynomial modified by a real parameter $\lambda$ and whose normalized version, $W_n^\lambda(z;a,b,c,d)$, satisfies the following recursion relation

$$\begin{aligned} zW_n^\lambda = & \left\{\frac{(n+a+b)(n+a+c)(n+a+d)(n+a+b+c+d-1)}{(2n+a+b+c+d)(2n+a+b+c+d-1)} + \frac{n(n+b+c-1)(n+b+d-1)(n+c+d-1)}{(2n+a+b+c+d-1)(2n+a+b+c+d-2)} - a^2 - \lambda\left[(n+a+c)(n+b+d)-n\right]\right\}W_n^\lambda \\ & + \frac{1}{2n+a+b+c+d-2}\sqrt{\frac{n(n+a+b-1)(n+c+d-1)(n+a+c-1)(n+a+d-1)(n+b+c-1)(n+b+d-1)(n+a+b+c+d-2)}{(2n+a+b+c+d-3)(2n+a+b+c+d-1)}}\,W_{n-1}^\lambda \\ & + \frac{1}{2n+a+b+c+d}\sqrt{\frac{(n+1)(n+a+b)(n+c+d)(n+a+c)(n+a+d)(n+b+c)(n+b+d)(n+a+b+c+d-1)}{(2n+a+b+c+d-1)(2n+a+b+c+d+1)}}\,W_{n+1}^\lambda \end{aligned} \tag{16}$$

which is the same as the original recursion relation (8) except for the presence of the new term, $-\lambda\left[(n+a+c)(n+b+d)-n\right]$, in the diagonal and sign reversal of the off-diagonal. Comparing (16) with the S3TRR (B10b) gives $f_n(z)=W_n^\lambda(z;\sigma-\tau,\sigma+\tau,\gamma,\gamma)$ and results in the following parameters assignments

$$\lambda=d\,,\quad 2\sigma=\mu+1\,,\quad 2\gamma=\nu+1\,,\quad 2\tau=\sqrt{(a+b+c-1)^2-4D}\,, \tag{17a}$$

$$z=-\frac{1}{4}(\nu+1)^2+\frac{\nu+1-ac}{2}-R-d\left[1+D+\frac{\mu+\nu}{2}-\frac{c}{2}(a+b)\right]. \tag{17b}$$

The nature of the solution depends on the value of $\lambda$ and on whether $D$ is greater than or less than $\frac{1}{4}(a+b+c-1)^2$. We do not have the expertise and needed mathematical skills to obtain the analytic properties of these modified Wilson polynomials. Nonetheless, numerical analysis shows that if $0<\lambda<1$ then the spectrum is pure continuous for all values of $D$ and $z\in\mathbb{R}$. The solution becomes

$$y(z,x)=p(\lambda,z)\sum_{n=0}^{\infty}W_n^\lambda(z;\sigma-\tau,\sigma+\tau,\gamma,\gamma)\phi_n(x) \tag{18}$$

where $p^2(\lambda,z)$ is the weight function of this modified Wilson polynomial,

$$\int_{-\infty}^{+\infty}p^2(\lambda,z)W_n^\lambda(z;a,b,c,d)W_n^\lambda(z;a,b,c,d)\,dz=\delta_{n,m}\,, \tag{19}$$

which is yet to be derived. On the other hand, if $\lambda\notin[0,1]$ and $4D>(a+b+c-1)^2$ then $\tau$ becomes pure imaginary and the spectrum is still pure continuous but with $z>0$ for $\lambda<0$ and $z<0$ for $\lambda>1$. However, if $4D<(a+b+c-1)^2$ then $\tau$ is real and the spectrum is a mix of a continuous part and a discrete part of finite size. For $\lambda<0$ (or $\lambda>1$) the discrete spectrum is negative (or positive) whereas the continuous spectrum is positive (or negative), respectively. Table 3 gives a summary of these scenarios. The solution for the mix spectrum case reads as follows

$$\begin{aligned} y_k(z,x) = p(\lambda,z)\sum\nolimits_{n=0}^{\infty} W_n^{\lambda}(z;\sigma-\tau,\sigma+\tau,\gamma,\gamma)\phi_n(x) \\ +q(\lambda,k)\sum\nolimits_{n=0}^{N} W_n^{\lambda}(z_k;\sigma-\tau,\sigma+\tau,\gamma,\gamma)\phi_n(x) \end{aligned} \tag{20}$$

where $q(\lambda,k)$ is obtained from the generalized orthogonality of $W_n^{\lambda}(z;a,b,c,d)$ that consists of a continuous integral and a discrete finite sum as follows

$$\int p^2(\lambda,z) W_n^{\lambda}(z;a,b,c,d) W_n^{\lambda}(z;a,b,c,d)\,dz + \sum_{k=0}^{N} q^2(\lambda,k) W_n^{\lambda}(z_k;a,b,c,d) W_n^{\lambda}(z_k;a,b,c,d) = \delta_{n,m} \tag{21}$$

We call upon experts in the field of orthogonal polynomials to derive the analytic properties of this modified version of the Wilson polynomial for all values of the deformation parameter $\lambda$.

## IV. The special solution

For this class of solutions, the differential equation parameters $\{A,B,C\}$ and the basis parameters $\{\alpha,\beta,\gamma,\mu,\nu\}$ are to be chosen and should satisfy the constraints given in the column headed "Case (a)" in Table 1.

We start by introducing a new four-parameter orthogonal polynomial $Q_n^{(\mu,\nu)}(z;\theta,\sigma)$ whose normalized version satisfies the following recursion relation

$$\begin{aligned} (\cos\theta) Q_n^{(\mu,\nu)}(z;\theta,\sigma) = \left\{ z\sin\theta \left[ \left(n+\tfrac{\mu+\nu+1}{2}\right)^2 - \sigma^2 \right]^{-1} + F_n \right\} Q_n^{(\mu,\nu)}(z;\theta,\sigma) \\ +2G_{n-1} Q_{n-1}^{(\mu,\nu)}(z;\theta,\sigma) + 2G_n Q_{n+1}^{(\mu,\nu)}(z;\theta,\sigma) \end{aligned} \tag{22}$$

where $z\in\mathbb{R}$, $0\le\theta<\pi$ and $\{F_n,G_n\}$ are defined in Appendix B by Eq. (B8).. This polynomial has a purely continuous spectrum and belongs to a new class of orthogonal polynomials introduced in Refs. [17,18]. One of the elements of this class, which was encountered in the past, satisfies the following recurrence relation

$$\begin{aligned} (\cos\theta) H_n^{(\mu,\nu)}(z;\theta,\sigma) = \left\{ z\sin\theta \left[ \left(n+\tfrac{\mu+\nu+1}{2}\right)^2 - \sigma^2 \right] + F_n \right\} H_n^{(\mu,\nu)}(z;\theta,\sigma) \\ +2G_{n-1} H_{n-1}^{(\mu,\nu)}(z;\theta,\sigma) + 2G_n H_{n+1}^{(\mu,\nu)}(z;\theta,\sigma) \end{aligned} \tag{23}$$

Taking $z=0$ turns (22) and (23) into the S3TRR of the (normalized) Jacobi polynomial $P_n^{(\mu,\nu)}(\cos\theta)$. If we compare (22) to the S3TRR (B10a) we obtain $f_n(z)=Q_n^{(\mu,\nu)}(z;\theta,\sigma)$ and the following parameters assignments

$$\cos\theta = 2d-1,\quad \sigma^2 = \tfrac{1}{4}(a+b+c-1)^2 - D,\quad z = \frac{R + dD - (c/2)\left[d(a+b)-a\right]}{\sqrt{d-d^2}}. \tag{24}$$

with the parameter constraint $0 \leq d \leq 1$. Therefore, the solution of the Heun-type differential equation (2) is written for a given value of *z* as

$$y(z,x) = p(z)\sum_{n=0}^{\infty} Q_n^{(\mu,\nu)}(z;\theta,\sigma)\phi_n(x), \tag{25}$$

where $\phi_n(x)$ is given by (5) and *p*(*z*) is a normalization factor, which is proportional to the square root of weight function of $Q_n^{(\mu,\nu)}(z;\theta,\sigma)$ that is yet to be derived.

If, on the other hand, $d \notin [0,1]$ then the spectrum depends on whether $\sigma$ is real or pure imaginary. For the latter case, the spectrum is continuous with $z > 0$. However, if $\sigma$ is real then the spectrum is a mix of a continuous positive spectrum and a discrete negative spectrum of finite size. We refer to this new polynomial by $G_n^{(\mu,\nu)}(k;\xi,\sigma)$ and its recursion relation is obtained from (22) by the replacement $\theta \to i\theta$ and $z \to -iz$ giving

$$\begin{aligned}\left(1+\xi^2\right)G_n^{(\mu,\nu)}(z;\xi,\sigma) = \left\{ z\left(1-\xi^2\right)\left[\left(n+\tfrac{\mu+\nu+1}{2}\right)^2 - \sigma^2\right]^{-1} + 2\xi F_n \right\} G_n^{(\mu,\nu)}(z;\xi,\sigma) \\ +4\xi\left[G_{n-1}G_{n-1}^{(\mu,\nu)}(z;\xi,\sigma) + G_n G_{n+1}^{(\mu,\nu)}(z;\xi,\sigma)\right]\end{aligned} \tag{26}$$

where $1 > \xi = e^{-\theta} > 0$. For real $\sigma$, the finite discrete spectrum $\{z_k\}_{k=0}^{N}$ is determined from the condition that forces the asymptotics ($n \to \infty$) of $Q_n^{(\mu,\nu)}(z;\theta,\sigma)$ to vanish (i.e., its spectrum formula). However, this asymptotics is yet to be derived analytically. Had this been known, we would have been able to write the solution of the Heun-type differential equation (2) for real $\sigma$ and a fixed value of the non-negative integer *k* as

$$y_k(z,x) = p(z)\sum_n G_n^{(\mu,\nu)}(z;\xi,\sigma)\phi_n(x) + q(k)\sum_{n=0}^{N} G_n^{(\mu,\nu)}(z_k;\xi,\sigma)\phi_n(x), \tag{27}$$

where *p*(*z*) and *q*(*k*) are obtained from the generalized orthogonality relation

$$\int_0^\infty p^2(z)\,G_n^{(\mu,\nu)}(z;\xi,\sigma)\,G_m^{(\mu,\nu)}(z;\xi,\sigma)\,dz + \sum_{k=0}^{N} q^2(k)\,G_n^{(\mu,\nu)}(z_k;\xi,\sigma)\,G_m^{(\mu,\nu)}(z_k;\xi,\sigma) = \delta_{n,m}. \tag{28}$$

# V. Conclusion

We introduced a nine-parameter ordinary second order linear differential equations of the Heun-type and derived three classes of its solution as infinite series of square integrable functions written in terms of the Jacobi polynomial. The expansion coefficients of the series are orthogonal polynomials satisfying S3TRR's. Some of these polynomials are new and were not treated previously in the mathematics literature while others are well-known hypergeometric type orthogonal polynomials (the Wilson polynomial or one of its modified versions).

The details of the "generalized" and "special" solution of the Heun-type differential equation (2) could not be given due to the lack of knowledge of the analytic properties of the corresponding orthogonal polynomials $W_n^\lambda(z;a,b,c,d)$ and $Q_n^{(\mu,\nu)}(z;\theta,\sigma)$, which are defined up to now by their S3TRR. Due to the prime significance of these polynomials in physics and mathematics, we call upon experts in the field to derive their properties (weight function, generating function, orthogonality, zeros, asymptotics, Rodrigues-type formula, differential or shift formula, etc.) and find their discrete versions.

# Appendix A: The Jacobi polynomial

For ease of reference, we list below the basic properties of the Jacobi polynomial, which is defined here by the replacement $x \to 2x-1$ in the classic definition as follows

$$P_n^{(\mu,\nu)}(x) = \frac{\Gamma(n+\mu+1)}{\Gamma(n+1)\Gamma(\mu+1)}\, {}_2F_1\left(-n, n+\mu+\nu+1; \mu+1; 1-x\right). \qquad \text{(A1)}$$

where $\mu > -1$ and $\nu > -1$ for $x \in [0,+1]$. It satisfies the following differential equation

$$\left\{x(1-x)\frac{d^2}{dx^2} + \left[\nu+1-x(\mu+\nu+2)\right]\frac{d}{dx} + n(n+\mu+\nu+1)\right\}P_n^{(\mu,\nu)}(x) = 0\,, \qquad \text{(A2)}$$

and differential relation

$$\begin{aligned} x(1-x)\frac{d}{dx}P_n^{(\mu,\nu)}(x) &= (n+\mu+\nu+1)\left[\frac{(\mu-\nu)n}{(2n+\mu+\nu)(2n+\mu+\nu+2)}P_n^{(\mu,\nu)}(x)\right. \\ &\left. + \frac{(n+\mu)(n+\nu)}{(2n+\mu+\nu)(2n+\mu+\nu+1)}P_{n-1}^{(\mu,\nu)}(x) - \frac{n(n+1)}{(2n+\mu+\nu+1)(2n+\mu+\nu+2)}P_{n+1}^{(\mu,\nu)}(x)\right] \end{aligned} \qquad \text{(A3)}$$

It also satisfies the following three-term recursion relation

$$\begin{aligned} x\,P_n^{(\mu,\nu)}(x) &= \frac{1}{2}\left[\frac{\nu^2-\mu^2}{(2n+\mu+\nu)(2n+\mu+\nu+2)} + 1\right]P_n^{(\mu,\nu)}(x) \\ &+ \frac{(n+\mu)(n+\nu)}{(2n+\mu+\nu)(2n+\mu+\nu+1)}P_{n-1}^{(\mu,\nu)}(x) + \frac{(n+1)(n+\mu+\nu+1)}{(2n+\mu+\nu+1)(2n+\mu+\nu+2)}P_{n+1}^{(\mu,\nu)}(x) \end{aligned} \qquad \text{(A4)}$$

The associated orthogonality relation reads as follows

$$\int_0^{+1} x^{\nu}(1-x)^{\mu}P_n^{(\mu,\nu)}(x)P_m^{(\mu,\nu)}(x)dx = \frac{1}{2n+\mu+\nu+1}\frac{\Gamma(n+\mu+1)\Gamma(n+\nu+1)}{\Gamma(n+1)\Gamma(n+\mu+\nu+1)}\delta_{nm} \qquad \text{(A5)}$$

# Appendix B: Deriving the S3TRR using the TRA

The requirement that the substitution of $y(x)$ given by (4) in the differential equations (2) result in a S3TRR for the expansion coefficients $\{f_n\}$ means that the action of the differential operator (2) on the basis element (5) must produce three terms containing $P_n^{(\mu,\nu)}(x)$ and $P_{n\pm1}^{(\mu,\nu)}(x)$ with $x$-independent multiplicative factors. Now, this differential operator action on $\phi_n(x)$ becomes

$$\mathcal{J}\phi_n(x) = c_n x^{\alpha}(x-1)^{\beta}(x-d)^{\gamma+1}$$

$$\left[x(x-1)\left(\frac{d}{dx}+\frac{\alpha}{x}+\frac{\beta}{x-1}+\frac{\gamma}{x-d}\right)^2 + x(x-1)\left(\frac{a}{x}+\frac{b}{x-1}+\frac{c}{x-d}\right)\times \right. \tag{B1}$$

$$\left.\left(\frac{d}{dx}+\frac{\alpha}{x}+\frac{\beta}{x-1}+\frac{\gamma}{x-d}\right)+\frac{1}{x-d}\left(\frac{A}{x}+\frac{B}{x-1}+\frac{C}{x-d}+xD-E\right)\right]P_n^{(\mu,\nu)}(x)$$

where $\mathcal{J} = x(x-1)(x-d)\mathcal{D}$ and the Heun-type equation (2) is written as $\mathcal{D}y(x) = 0$. Expanding and collecting similar terms give

$$\mathcal{J}\phi_n(x) = c_n x^{\alpha}(x-1)^{\beta}(x-d)^{\gamma+1}\left\{x(x-1)\left[\frac{d^2}{dx^2}+\left(\frac{2\alpha+a}{x}+\frac{2\beta+b}{x-1}+\frac{2\gamma+c}{x-d}\right)\frac{d}{dx}\right.\right.$$

$$\left.+\frac{\alpha(\alpha-1+a)}{x^2}+\frac{\beta(\beta-1+b)}{(x-1)^2}+\frac{\gamma(\gamma-1+c)}{(x-d)^2}+\frac{2\alpha\beta+\alpha b+\beta a}{x(x-1)}+\frac{2\alpha\gamma+\alpha c+\gamma a}{x(x-d)}+\frac{2\beta\gamma+\beta c+\gamma b}{(x-1)(x-d)}\right] \tag{B2}$$

$$\left.+\frac{1}{x-d}\left(\frac{A}{x}+\frac{B}{x-1}+\frac{C}{x-d}+xD-E\right)\right\}P_n^{(\mu,\nu)}(x)$$

Using the differential equation of the Jacobi polynomials (A2) turns this into the following

$$\mathcal{J}\phi_n(x) = c_n x^{\alpha}(x-1)^{\beta}(x-d)^{\gamma+1}\left\{\left(\frac{2\alpha+a-\nu-1}{x}+\frac{2\beta+b-\mu-1}{x-1}+\frac{2\gamma+c}{x-d}\right)x(x-1)\frac{d}{dx}\right.$$

$$+n\left(n+\mu+\nu+1\right)+\left(\alpha+\beta+\gamma\right)\left(\alpha+\beta+\gamma+a+b+c-1\right)-\frac{\alpha(\alpha-1+a)}{x}+\frac{\beta(\beta-1+b)}{x-1}+d(d-1)\frac{\gamma(\gamma-1+c)}{(x-d)^2} \tag{B3}$$

$$\left.+\frac{1}{x-d}\left[(2d-1)\gamma\left(\gamma-1+c\right)+(d-1)\left(2\alpha\gamma+\alpha c+\gamma a\right)+d\left(2\beta\gamma+\beta c+\gamma b\right)+\frac{A}{x}+\frac{B}{x-1}+\frac{C}{x-d}+xD-E\right]\right\}P_n^{(\mu,\nu)}(x)$$

where we have also used the identities

$$\frac{x}{x-d}=1+\frac{d}{x-d},\quad \frac{x-1}{x-d}=1+\frac{d-1}{x-d},\quad \frac{x(x-1)}{(x-d)^2}=1+\frac{2d-1}{x-d}+\frac{d(d-1)}{(x-d)^2}. \tag{B4}$$

The differential property of the Jacobi polynomials (A3) turns (B3) into the following

$$\mathcal{J}\phi_n(x) = c_n x^{\alpha}(x-1)^{\beta}(x-d)^{\gamma+1}\left\{\left(n+\mu+\nu+1\right)\left(\frac{2\alpha+a-\nu-1}{x}+\frac{2\beta+b-\mu-1}{x-1}+\frac{2\gamma+c}{x-d}\right)\times\right.$$

$$\left[\frac{(\nu-\mu)n}{(2n+\mu+\nu)(2n+\mu+\nu+2)}P_n^{(\mu,\nu)}-\frac{(n+\mu)(n+\nu)}{(2n+\mu+\nu)(2n+\mu+\nu+1)}P_{n-1}^{(\mu,\nu)}+\frac{n(n+1)}{(2n+\mu+\nu+1)(2n+\mu+\nu+2)}P_{n+1}^{(\mu,\nu)}\right]$$

$$+\left[n\left(n+\mu+\nu+1\right)+\left(\alpha+\beta+\gamma\right)\left(\alpha+\beta+\gamma+a+b+c-1\right)\right]P_n^{(\mu,\nu)} \tag{B5}$$

$$+\frac{1}{x-d}\left[-\alpha\left(\alpha-1+a\right)+\beta\left(\beta-1+b\right)+(2d-1)\gamma\left(\gamma-1+c\right)+(d-1)\left(2\alpha\gamma+\alpha c+\gamma a\right)+d\left(2\beta\gamma+\beta c+\gamma b\right)\right]P_n^{(\mu,\nu)}$$

$$\left.+\frac{1}{x-d}\left[\alpha(\alpha-1+a)\frac{d}{x}+\beta(\beta-1+b)\frac{1-d}{x-1}+d(d-1)\frac{\gamma(\gamma-1+c)}{x-d}+\frac{A}{x}+\frac{B}{x-1}+\frac{C}{x-d}+xD-E\right]P_n^{(\mu,\nu)}\right\}$$

The tridiagonal requirement and recursion relation of the Jacobi polynomials (A4) allow only $x$-independent terms to multiply $P_{n\pm1}^{(\mu,\nu)}(x)$ in (B5) and only linear terms in $x$ multiplying $P_n^{(\mu,\nu)}(x)$. These constraints lead to the four set of scenarios shown in Table 1. The two cases corresponding to the first and second columns of Table 1 are obtained by multiplying the content of the curly brackets in (B5) by $x-d$ then imposing the tridiagonal requirement. That is, factoring out $\frac{1}{x-d}$ making the overall factor multiplying (B5) equal to $c_n x^{\alpha}(x-1)^{\beta}(x-d)^{\gamma}$. The case in the third column of Table 1 is obtained by multiplying the content of the curly brackets in (B5) by $x-1$ then imposing the tridiagonal requirement. That is, factoring out $\frac{1}{x-1}$ making the overall multiplicative factor $c_n x^{\alpha}(x-1)^{\beta-1}(x-d)^{\gamma+1}$. The fourth column of the Table, on the other hand, is obtained by multiplying the content of the curly brackets by *x* then imposing the tridiagonal requirement. That is, we factor out $\frac{1}{x}$ making the overall multiplicative factor $c_n x^{\alpha-1}(x-1)^{\beta}\ (x-d)^{\gamma+1}$. The top five rows of the Table relate the basis parameters $\{\alpha,\beta,\gamma,\mu,\nu\}$ to the differential equation parameters $\{a,b,c,d,A,B,C\}$. The bottom five rows are constraints on, or relations among, the differential equation parameters brought about either by reality, which results in the inequalities, or by the tridiagonal requirement resulting in the equalities. Substitution of the parameters relations for each of the four cases from the Table into Eq. (B5) maps it into the following corresponding relations

$$\mathcal{J}\phi_n(x) = c_n x^{\alpha}(x-1)^{\beta}(x-d)^{\gamma}\left\{\tfrac{1}{4}\left[\left(2n+\mu+\nu+1\right)^2-\left(a+b+c-1\right)^2\right](x-d)+xD \right.$$
$$\left. +\tfrac{A}{d}-\tfrac{B}{1-d}+\tfrac{2d-1}{1-d}\tfrac{C}{d}-E-\tfrac{c}{2}\left[d\left(a+b\right)-a\right]\right\}P_n^{(\mu,\nu)} \tag{B6a}$$

$$\mathcal{J}\phi_n(x) = c_n x^{\alpha}(x-1)^{\beta}(x-d)^{\gamma}\left\{g\left(n+\mu+\nu+1\right)\times\right.$$
$$\left[\frac{(\nu-\mu)n}{(2n+\mu+\nu)(2n+\mu+\nu+2)}P_n^{(\mu,\nu)}-\frac{(n+\mu)(n+\nu)}{(2n+\mu+\nu)(2n+\mu+\nu+1)}P_{n-1}^{(\mu,\nu)}+\frac{n(n+1)}{(2n+\mu+\nu+1)(2n+\mu+\nu+2)}P_{n+1}^{(\mu,\nu)}\right]$$
$$+\tfrac{1}{4}\left[\left(2n+\mu+\nu+1\right)^2-\left(a+b+c-1\right)^2+g^2+2g\left(\mu+\nu+1\right)\right](x-d)P_n^{(\mu,\nu)}+xDP_n^{(\mu,\nu)}$$
$$\left.+\left(\tfrac{A}{d}-\tfrac{B}{1-d}+\tfrac{2d-1}{1-d}\tfrac{C}{d}-E-\tfrac{c}{2}\left[d\left(a+b\right)-a\right]+\tfrac{g}{2}\left[d\left(\mu+\nu+2\right)-\nu-1\right]\right)P_n^{(\mu,\nu)}\right\} \tag{B6b}$$

where $g=2\gamma+c$.

$$\mathcal{J}\phi_n(x) = c_n x^{\alpha}(x-1)^{\beta-1}(x-d)^{\gamma+1}\left\{p\left(n+\mu+\nu+1\right)\times\right.$$
$$\left[\frac{(\nu-\mu)n}{(2n+\mu+\nu)(2n+\mu+\nu+2)}P_n^{(\mu,\nu)}-\frac{(n+\mu)(n+\nu)}{(2n+\mu+\nu)(2n+\mu+\nu+1)}P_{n-1}^{(\mu,\nu)}+\frac{n(n+1)}{(2n+\mu+\nu+1)(2n+\mu+\nu+2)}P_{n+1}^{(\mu,\nu)}\right]$$
$$\left.+\tfrac{1}{4}\left[\left(2n+\mu+\nu+1\right)^2-\left(a+b+c-1\right)^2+p^2+2p\left(\mu+\nu+1\right)\right](x-1)P_n^{(\mu,\nu)}\right\} \tag{B6c}$$

where $p=2\beta+b-\mu-1$.

$$\mathcal{J}\phi_n(x) = c_n x^{\alpha-1}(x-1)^{\beta}(x-d)^{\gamma+1}\left\{q\left(n+\mu+\nu+1\right)\times\right.$$
$$\left[\frac{(\nu-\mu)n}{(2n+\mu+\nu)(2n+\mu+\nu+2)}P_n^{(\mu,\nu)}-\frac{(n+\mu)(n+\nu)}{(2n+\mu+\nu)(2n+\mu+\nu+1)}P_{n-1}^{(\mu,\nu)}+\frac{n(n+1)}{(2n+\mu+\nu+1)(2n+\mu+\nu+2)}P_{n+1}^{(\mu,\nu)}\right]$$
$$\left.+\tfrac{1}{4}\left[\left(2n+\mu+\nu+1\right)^2-\left(a+b+c-1\right)^2+q^2+2q\left(\mu+\nu+1\right)\right]xP_n^{(\mu,\nu)}\right\} \tag{B6d}$$

where $q=2\alpha+a-\nu-1$. Next, we use the three-term recursion relation of the Jacobi polynomials (A4) in the above four relations and write them in terms of the basis elements $\{\phi_n(x)\}$ instead of $P_n^{(\mu,\nu)}(x)$ (i.e., taking care of the normalization factor $c_n$). Requiring symmetry of the tridiagonal representation (i.e., the *n*-dependent factor multiplying $\phi_{n+1}(x)$ must be identical to that multiplying $\phi_{n-1}(x)$ but with $n\to n+1$) dictates that $g=p=q=1$ and gives

$$\mathcal{J}\phi_n(x)=\Bigg\{R-\frac{c}{2}\big[d(a+b)-a\big]-d\left[\left(n+\tfrac{\mu+\nu+1}{2}\right)^2-\frac{1}{4}(a+b+c-1)^2\right]$$
$$+\frac{1}{2}(F_n+1)\left[\left(n+\tfrac{\mu+\nu+1}{2}\right)^2-\frac{1}{4}(a+b+c-1)^2+D\right]\Bigg\}\phi_n(x) \tag{B7a}$$
$$+G_{n-1}\left[\left(n+\tfrac{\mu+\nu+1}{2}\right)^2-\frac{1}{4}(a+b+c-1)^2+D\right]\phi_{n-1}(x)+G_n\left[\left(n+\tfrac{\mu+\nu+1}{2}\right)^2-\frac{1}{4}(a+b+c-1)^2+D\right]\phi_{n+1}(x)$$

$$\mathcal{J}\phi_n(x)=\Bigg\{-\frac{n(n+\mu)}{2n+\mu+\nu}+\frac{1}{2}(F_n+1)\left[\left(n+\tfrac{\mu+\nu}{2}+1\right)^2-\frac{1}{4}(a+b+c-1)^2+D\right]$$
$$-d\left[\left(n+\tfrac{\mu+\nu+1}{2}\right)^2-\frac{1}{4}(a+b+c-1)^2+\frac{c}{2}(a+b)-\frac{1}{4}\right]-\frac{\nu+1}{2}+\frac{ac}{2}+R\Bigg\}\phi_n(x) \tag{B7b}$$
$$+G_{n-1}\left[\left(n+\tfrac{\mu+\nu}{2}\right)^2-\frac{1}{4}(a+b+c-1)^2+D\right]\phi_{n-1}(x)+G_n\left[\left(n+\tfrac{\mu+\nu}{2}+1\right)^2-\frac{1}{4}(a+b+c-1)^2+D\right]\phi_{n+1}(x)$$

$$\frac{x-1}{x-d}\mathcal{J}\phi_n(x)=\left\{\frac{n(n+\nu)}{2n+\mu+\nu}+\frac{1}{2}(F_n-1)\left[\left(n+\tfrac{\mu+\nu}{2}+1\right)^2-\frac{1}{4}(a+b+c-1)^2\right]\right\}\phi_n(x)$$
$$+G_{n-1}\left[\left(n+\tfrac{\mu+\nu}{2}\right)^2-\frac{1}{4}(a+b+c-1)^2\right]\phi_{n-1}(x)+G_n\left[\left(n+\tfrac{\mu+\nu}{2}+1\right)^2-\frac{1}{4}(a+b+c-1)^2\right]\phi_{n+1}(x) \tag{B7c}$$

$$\frac{x}{x-d}\mathcal{J}\phi_n(x)=\left\{-\frac{n(n+\mu)}{2n+\mu+\nu}+\frac{1}{2}(F_n+1)\left[\left(n+\tfrac{\mu+\nu}{2}+1\right)^2-\frac{1}{4}(a+b+c-1)^2\right]\right\}\phi_n(x)$$
$$+G_{n-1}\left[\left(n+\tfrac{\mu+\nu}{2}\right)^2-\frac{1}{4}(a+b+c-1)^2\right]\phi_{n-1}(x)+G_n\left[\left(n+\tfrac{\mu+\nu}{2}+1\right)^2-\frac{1}{4}(a+b+c-1)^2\right]\phi_{n+1}(x) \tag{B7d}$$

where $R=\frac{A}{d}-\frac{B}{1-d}+\frac{2d-1}{1-d}\frac{C}{d}-E$ and

$$F_n=\frac{\nu^2-\mu^2}{(2n+\mu+\nu)(2n+\mu+\nu+2)},\quad G_n=\frac{1}{2n+\mu+\nu+2}\sqrt{\frac{(n+1)(n+\mu+1)(n+\nu+1)(n+\mu+\nu+1)}{(2n+\mu+\nu+1)(2n+\mu+\nu+3)}}\,. \tag{B8}$$

In this calculation, we employed the following identities

$$\frac{2n(n+\mu+\nu+1)(\mu-\nu)}{(2n+\mu+\nu)(2n+\mu+\nu+2)}=\frac{2n(n+\mu)}{2n+\mu+\nu}-n(F_n+1), \tag{B9b}$$

$$\frac{2n(n+\mu+\nu+1)(\nu-\mu)}{(2n+\mu+\nu)(2n+\mu+\nu+2)}=\frac{2n(n+\nu)}{2n+\mu+\nu}+n(F_n-1). \tag{B9c}$$

Putting all above findings together, we obtain the revised basis parameter assignments for cases (c) and (d) as shown in Table 2. Using these results and (B7a)-(B7d) in the differential equation $\mathcal{J}y(x)=0$ we obtain finally the following S3TRR for the expansion coefficients $\{f_n\}$

$$\left(2d-1\right)f_n(z)=\left\{\frac{2\left(dD+R\right)-c\left[d\left(a+b\right)-a\right]}{\left(n+\frac{\mu+\nu+1}{2}\right)^2-\frac{1}{4}\left(a+b+c-1\right)^2+D}+F_n\right\}f_n(z) \qquad \text{(B10a)}$$
$$+2G_{n-1}f_{n-1}(z)+2G_n f_{n+1}(z)$$

$$\left\{\frac{\nu+1-ac}{2}-R-d\left[1+D+\frac{\mu+\nu}{2}-\frac{c}{2}\left(a+b\right)\right]\right\}f_n(z)=$$
$$\left\{-\frac{n(n+\mu)}{2n+\mu+\nu}+dn+\frac{1}{2}\left(F_n+1-2d\right)\left[\left(n+\tfrac{\mu+\nu}{2}+1\right)^2-\frac{1}{4}\left(a+b+c-1\right)^2+D\right]\right\}f_n(z) \qquad \text{(B10b)}$$
$$+G_{n-1}\left[\left(n+\tfrac{\mu+\nu}{2}\right)^2-\frac{1}{4}\left(a+b+c-1\right)^2+D\right]f_{n-1}(z)+G_n\left[\left(n+\tfrac{\mu+\nu}{2}+1\right)^2-\frac{1}{4}\left(a+b+c-1\right)^2+D\right]f_{n+1}(z)$$

$$0=\left\{\frac{n\left(n+\nu\right)}{2n+\mu+\nu}+\frac{1}{2}\left(F_n-1\right)\left[\left(n+\tfrac{\mu+\nu}{2}+1\right)^2-\frac{1}{4}\left(a+b+c-1\right)^2\right]\right\}f_n(z)$$
$$+G_{n-1}\left[\left(n+\tfrac{\mu+\nu}{2}\right)^2-\frac{1}{4}\left(a+b+c-1\right)^2\right]f_{n-1}(z)+G_n\left[\left(n+\tfrac{\mu+\nu}{2}+1\right)^2-\frac{1}{4}\left(a+b+c-1\right)^2\right]f_{n+1}(z) \qquad \text{(B10c)}$$

$$0=\left\{\frac{n\left(n+\mu\right)}{2n+\mu+\nu}-\frac{1}{2}\left(F_n+1\right)\left[\left(n+\tfrac{\mu+\nu}{2}+1\right)^2-\frac{1}{4}\left(a+b+c-1\right)^2\right]\right\}f_n(z)$$
$$-G_{n-1}\left[\left(n+\tfrac{\mu+\nu}{2}\right)^2-\frac{1}{4}\left(a+b+c-1\right)^2\right]f_{n-1}(z)-G_n\left[\left(n+\tfrac{\mu+\nu}{2}+1\right)^2-\frac{1}{4}\left(a+b+c-1\right)^2\right]f_{n+1}(z) \qquad \text{(B10d)}$$

Taking $f_0(z)=1$ these S3TRR's give $f_n(z)$ as a polynomial of degree *n* in some proper argument *z*, which is a function of the differential equation parameters. The solutions of these S3TRR's are obtained in sections II-IV in terms of orthogonal polynomials with continuous and/or discrete spectra.

It is worth noting that the S3TRR (B10d) is obtained from (B10c) by the following parameter exchange/map

$$\mu\leftrightarrow\nu,\quad \alpha\leftrightarrow\beta,\quad a\leftrightarrow b,\quad \frac{A}{d}\leftrightarrow\frac{B}{1-d},\quad E\to E-\frac{2A}{d}+\frac{2B}{1-d}+\frac{c}{2}\left(b-a\right). \qquad \text{(B11)}$$

in addition to $f_n\to(-1)^n f_n$, which is equivalent to $x\to 1-x$. Consequently, in the solution of the Heun-type differential equation we consider the S3TRR (B10a), (B10b) and either (B10c) or (B10d).

# Tables Captions:

**Table 1**: The top five rows give the values of the basis parameters $\{\alpha,\beta,\gamma,\mu,\nu\}$ in terms of the Heun-type equation parameters $\{a,b,c,d,A,B,C\}$ for each of the four classes of its solution. The bottom five rows give constraints on the equation parameters for each solution class.

**Table 2**: The revised values of the basis parameters $\{\alpha,\beta,\gamma,\mu,\nu\}$ from Table 1 for the two restricted classes of solutions.

**Table 3**: The type of spectrum of the generalized solution obtained in Sec. III as a function of the deformation parameter $\lambda$ of the Wilson polynomial, $W_n^{\lambda}(z;\sigma-\tau,\sigma+\tau,\gamma,\gamma)$, and reality of the parameter $\tau$.

**Table 2**

| **Case (c)**<br>**Restricted Solution** | **Case (d)**<br>**Restricted Solution** |
|---|---|
| $2\alpha = \nu + 1 - a$ | $2\alpha = \nu + 2 - a$ |
| $2\beta = \mu + 2 - b$ | $2\beta = \mu + 1 - b$ |
| $2\gamma = -c$ | $2\gamma = -c$ |
| $\nu^2 = (1-a)^2 - 4\dfrac{A}{d}$ | $(\nu+1)^2 = (1-a)^2 - 4\dfrac{A}{d}$ |
| $(\mu+1)^2 = (1-b)^2 - \dfrac{4B}{1-d}$ | $\mu^2 = (1-b)^2 - \dfrac{4B}{1-d}$ |

**Table 3**

| | **Re(*τ*) = 0** | **Im(*τ*) = 0** |
|---|---|---|
| $0 < \lambda < 1$ | $z \in \mathbb{R}$ | $z \in \mathbb{R}$ |
| $\lambda < 0$ | $z > 0$ | $z > 0,\ z_k < 0$ |
| $\lambda > 1$ | $z < 0$ | $z < 0,\ z_k > 0$ |

**Table 1**

| **Case (a)**<br>**Special Solution** | **Case (b)**<br>**Generalized Solution** | **Case (c)**<br>**Restricted Solution** | **Case (d)**<br>**Restricted Solution** |
|---|---|---|---|
| $2\alpha = \nu + 1 - a$ | $2\alpha = \nu + 1 - a$ | $2\alpha = \nu + 1 - a$ | $\left(2\alpha + a - 1\right)^2 = (1-a)^2 - 4\frac{A}{d}$ |
| $2\beta = \mu + 1 - b$ | $2\beta = \mu + 1 - b$ | $\left(2\beta + b - 1\right)^2 = (1-b)^2 - \frac{4B}{1-d}$ | $2\beta = \mu + 1 - b$ |
| $2\gamma = -c$ | $\left(2\gamma + c - 1\right)^2 = (1-c)^2 + \frac{4C/d}{1-d}$ | $2\gamma = -c$ | $2\gamma = -c$ |
| $\nu^2 = (1-a)^2 - 4\frac{A}{d}$ | $\nu^2 = (1-a)^2 - 4\frac{A}{d}$ | $\nu^2 = (1-a)^2 - 4\frac{A}{d}$ | ××××× |
| $\mu^2 = (1-b)^2 - \frac{4B}{1-d}$ | $\mu^2 = (1-b)^2 - \frac{4B}{1-d}$ | ××××× | $\mu^2 = (1-b)^2 - \frac{4B}{1-d}$ |
| ××××× | ××××× | $D = 0$ | $D = 0$ |
| $1 = (1-c)^2 + \frac{4C/d}{1-d}$ | $\frac{4C/d}{1-d} \geq -(1-c)^2$ | $1 = (1-c)^2 + \frac{4C/d}{1-d}$ | $1 = (1-c)^2 + \frac{4C/d}{1-d}$ |
| ××××× | ××××× | $\frac{A}{d} - \frac{B}{1-d} + \frac{2d-1}{1-d}\frac{C}{d} - E = \frac{c}{2}\left[d\left(a+b\right) - a\right]$ | $\frac{A}{d} - \frac{B}{1-d} + \frac{2d-1}{1-d}\frac{C}{d} - E = \frac{c}{2}\left[d\left(a+b\right) - a\right]$ |
| $4\frac{A}{d} \leq (1-a)^2$ | $4\frac{A}{d} \leq (1-a)^2$ | $4\frac{A}{d} \leq (1-a)^2$ | $4\frac{A}{d} \leq (1-a)^2$ |
| $\frac{4B}{1-d} \leq (1-b)^2$ | $\frac{4B}{1-d} \leq (1-b)^2$ | $\frac{4B}{1-d} \leq (1-b)^2$ | $\frac{4B}{1-d} \leq (1-b)^2$ |